\documentclass[11pt,reqno]{amsart}
\usepackage{amssymb,amscd,amsbsy}
\usepackage{amssymb,amscd,amsbsy,mathrsfs}
\newcommand{\lb}{\linebreak}

\newcommand{\s}{\sigma}

\newcommand{\f}{\varphi}

\newcommand{\D}{\Delta}
\renewcommand{\L}{\Lambda}

\newcommand{\E}{{\mathscr E}}

\newcommand{\h}{{\mathscr H}}

\newcommand{\X}{{\mathscr X}}
\newcommand{\Y}{{\mathscr Y}}

\newcommand{\R}{{\Bbb R}}

\newcommand{\bs}{\boldsymbol}

\newcommand{\m}{{\boldsymbol m}}
\newcommand{\bS}{{\boldsymbol S}}

\newcommand{\rf}[1]{(\ref{#1})}

\newcommand{\df}{\stackrel{\mathrm{def}}{=}}

\newcommand{\trace}{\operatorname{trace}}

\newcommand{\eeq}{\end{equation}}
\newcommand{\beq}{\begin{equation}}
\newcommand{\bay}{\begin{eqnarray}}
\newcommand{\ba}{\begin{align*}}
\newcommand{\ea}{\end{align*}}
\newcommand{\ey}{\end{eqnarray}}
\newcommand{\bey}{\begin{eqnarray*}}
\newcommand{\eey}{\end{eqnarray*}}

\newcommand{\be}{\infty}

\newcommand{\bl}{\blacksquare}

\newcommand{\Pf}{{\bf Proof. }}

\newtheorem{thm}{\hspace{\parindent}Theorem}[section]

\theoremstyle{remark}

\newtheorem*{rem*}{Remark}

\newcommand{\OL}{{\rm OL}}

\newcommand\Li{{\rm Lip}}
\newcommand\fM{\frak M}

\newcommand\dg{\frak D}

\begin{document}

\newcommand{\vse}{\vspace{.2in}}
\numberwithin{equation}{section}

\title{The Lifshits--Krein trace formula\\ and operator Lipschitz functions}
\author{V.V. Peller}
\thanks{the author is partially supported by NSF grant DMS 1300924}

\begin{abstract}
We describe the maximal class of functions $f$ on the real line, for which 
the Lifshitz--Krein trace formula $\trace(f(A)-f(B))=\int_\R f'(s)\bs{\xi}(s)\,ds$
holds for arbitrary self-adjoint operators $A$ and $B$ with $A-B$ in the trace class $\bS_1$. We prove that this class of functions coincide with the class of operator Lipschitz functions. 
\end{abstract}

\maketitle

\tableofcontents
\normalsize

\

%
%
%
%

\setcounter{section}{0}
\section{\bf Introduction}
\setcounter{equation}{0}
\label{In}

\

The purpose of this paper is to describe the class functions, for which the Lifshitz--Krein trace formula holds. The Lifshitz--Krein trace formula plays a significant role in perturbation theory. It was discovered by Lifshits \cite{L} in a special case and by Krein \cite{Kr} in the general case. This formula allows one to compute the trace of the difference
$f(A)-f(B)$ of a function $f$ of an unperturbed self-adjoint operator $A$ and a perturbed self-adjoint operator $B$
provided the perturbation $B-A$ belongs to trace class $\bS_1$.
M.G. Krein proved that for each such pair there exists a unique function $\bs{\xi}$ in $L^1(\R)$ such that for every function $f$ whose derivative is
the Fourier transform of and $L^1$ function, the operator $f(A)-f(B)$ belongs to $\bS_1$ and the following trace formula holds:
\bay
\label{LKf}
\trace\big(f(A)-f(B)\big)=\int_\R f'(s)\bs{\xi}(s)\,ds
\ey
(see \cite{Kr}). The function $\bs{\xi}$ is called the {\it spectral shift function associated with the pair} $(A,B)$. Clearly, the right-hand side of \rf{LKf} makes sense for arbitrary Lipschitz function $f$. In this connection Krein asked the question of whether it is true that for an arbitrary Lipschitz function $f$, the operator $f(A)-f(B)$ is in $\bS_1$ and trace formula \rf{LKf} holds.
{\it It turns out that this is false}. In \cite{F2} Farforovskaya gave an example of self-adjoint operators $A$ and $B$ with $A-B\in\bS_1$ and a Lipschitz function $f$ such that $f(A)-f(B)\not\in\bS_1$. 

Later it was shown in \cite{Pe2} and \cite{Pe3} that formula \rf{LKf} holds, whenever $A$ and $B$ are self-adjoint operators with $A-B\in\bS_1$ and $f$ belongs to the Besov space $B_{\be,1}^1(\R)$ (we refer the reader to \cite{Pee} and \cite{Pe4} for an introduction to Besov classes). Necessary conditions are also obtained in \cite{Pe2} and \cite{Pe3}. In particular, it was shown in \cite{Pe2} and \cite{Pe3} that 
if $f(A)-f(B)\in\bS_1$ whenever $A$ and $B$ are self-adjoint operators with $A-B\in\bS_1$, then $f$ locally belongs to the Besov space $B_{1,1}^1(\R)$. Note that those necessary conditions were deduced from the description of trace class 
Hankel operators \cite{Pe1} (see also \cite{Pe4}).

The main objective of this paper is to describe the class of functions $f$, for which
trace formula \rf{LKf} holds for arbitrary self-adjoint operators $A$ and $B$ with $A-B\in\bS_1$.

It is well known (see e.g. \cite{AP}) that for a function $f$ on $\R$, the following properties are equivalent:

(i) {\it there exists a positive number $C$ such that 
\bay
\label{OLn}
\|f(A)-f(B)\|\le C\|A-B\|
\ey
for all bounded self-adjoint operators $A$ and $B$;}

(ii) {\it there exists a positive number $C$ such that inequality {\em\rf{OLn}} holds,
whenever $A$ and $B$ are (not necessarily bounded) self-adjoint operators such that $A-B$ is bounded;}

(iii) {\it there exists a positive number $C$ such that 
\bay
\label{OLyad}
\|f(A)-f(B)\|_{\bS_1}\le C\|A-B\|_{\bS_1}
\ey
for all bounded self-adjoint operators $A$ and $B$ with $A-B\in\bS_1$;}

(iv) {\it there exists a positive number $C$ such that inequality {\em\rf{OLyad}} holds,,
whenever $A$ and $B$ are (not necessarily bounded) self-adjoint operators such that 
$A-B\in\bS_1$;}

(v) $f(A)-f(B)\in\bS_1${\it, whenever $A$ and $B$ are (not necessarily bounded) self-adjoint operators such that $A-B\in\bS_1$.}

Note that the minimal value of the constant $C$ is the same in (i)--(iv).

Functions satisfying (i) are called {\it operator Lipschitz}. We denote by ${\rm OL}(\R)$ the space of operator Lipschitz functions on $\R$. For 
$f\in{\rm OL}(\R)$, we define its quasi norm $\|f\|_{\rm OL}$ as the infimum of all constants $C$, for which inequality \rf{OLn} holds.
In other words,
\begin{align*}
\|f\|_{\rm OL}&=
\sup\left\{\frac{\|f(A)-f(B)\|_{\bS_1}}{\|A-B\|_{\bS_1}}:
~A~\mbox{ and }~B~\mbox{ are self-adjoint},~A-B~\mbox{ is bounded}\right\}\\[.2cm]
&=\sup\left\{\frac{\|f(A)-f(B)\|}{\|A-B\|}:
~A~\mbox{ and }~B~\mbox{ are self-adjoint},~A-B\in\bS_1\right\}.
\end{align*}

It was shown in \cite{JW} that operator Lipschitz functions are differentiable everywhere on $\R$. Note that this implies that the function $x\mapsto|x|$
is not operator Lipschitz, the fact established earlier in \cite{Mc} and \cite{Ka}. On the other hand, an operator Lipschitz function does not have to be continuously differentiable; in particular, the function $x\mapsto x^2\sin x^{-1}$
is operator Lipschitz, see \cite{KS}.

For a differentiable function $f$ on $\R$, we consider the divided difference
$\dg f$ defined by
\bay
\label{razdraz}
(\dg f)(x,y)=\left\{\begin{array}{ll}\frac{f(x)-f(y)}{x-y},&x\ne y\\[.2cm]
f'(x),&x=y.\end{array}
\right.
\ey
It turns out (see e.g., \cite{AP}) that a differentiable function on $\R$ is operator Lipschitz if and only if the divided difference $\dg f$ is a Schur multiplier (see \S\,\ref{dois}) for the definition.

The main purpose of this paper is to prove that {\it the condition $f\in{\rm OL}(\R)$
is not only a necessary condition for the Lifshits--Krein trace formula {\em\rf{LKf}}
to hold for arbitrary self-adjoint operators $A$ and $B$ with $A-B\in\bS_1$, but is also sufficient.} This will be proved in \S\,\ref{self}.

In \S~\ref{dois} We define double operator integrals and Schur multipliers.
In \S~\ref{HSdif} we state a result of \cite{KPSS} on the differentiability of the function $t\mapsto f(A+tK)-f(A)$ in the Hilbert--Schmidt norm. We state a characterization of the space of Schur multipliers in terms of Haagerup tensor products in \S~\ref{Haag}.
Finally, in \S~\ref{sleddoi} we obtain a formula for the trace of double operator integrals.


%
%
%

\

\section{\bf Double operator integrals and Schur multipliers}
\setcounter{equation}{0}
\label{dois}

\

Double operator integrals appeared in the paper \cite{DK} by Daletskii and S.G. Krein. 
Later the beautiful theory of double operator integrals was created by Birman and Solomyak in \cite{BS1}, \cite{BS2}, and \cite{BS4}. 

Let $(\X,E_1)$ and $(\Y,E_2)$ be spaces with spectral measures $E_1$ and $E_2$
on a Hilbert space $\h$ and let $\Phi$ be a bounded measurable function on $\X\times\Y$. Double operator integrals are expressions of the form
\bay
\label{doi}
\int\limits_\X\int\limits_\Y\Phi(x,y)\,d E_1(x)T\,dE_2(y).
\ey
Birman and Solomyak starting point is the case when $T$ belongs to the Hilbert--Schmidt class $\bS_2$. For a bounded measurable function on $\Phi$ on $\X\times\Y$ and an operator $T$ of class $\bS_2$, consider the spectral measure $\E$ whose values are orthogonal
projections on the Hilbert space $\bS_2$, which is defined by
$$
\E(\L\times\D)T=E_1(\L)TE_2(\D),\quad T\in\bS_2,
$$
$\L$ and $\D$ being measurable subsets of $\X$ and $\Y$. 
It was shown in \cite{BS} that $\E$ extends to a spectral measure on
$\X\times\Y$. For a bounded measurable function $\Phi$ on $\X\times\Y$, the double operator integral \rf{doi} is defined by
$$
\int\limits_\X\int\limits_\Y\Phi(x,y)\,d E_1(x)T\,dE_2(y)\df
\left(\,\,\int\limits_{\X\times\Y}\Phi\,d\E\right)T.
$$
Clearly,
$$
\left\|\int\limits_\X\int\limits_\Y\Phi(x,y)\,dE_1(x)T\,dE_2(y)\right\|_{\bS_2}
\le\|\Phi\|_{L^\be}\|T\|_{\bS_2}.
$$
If
$$
\int\limits_\X\int\limits_\Y\Phi(x,y)\,d E_1(x)T\,dE_2(y)\in\bS_1
$$
for every $T\in\bS_1$, we say that $\Phi$ is a {\it Schur multiplier of $\bS_1$ associated with
the spectral measures $E_1$ and $E_2$}. {\it We denote by $\fM(E_1,E_2)$ the space of
Schur multipliers of $\bS_1$ with respect to $E_1$ and $E_2$.} The norm 
$\|\Phi\|_{\fM(E_1,E_2)}$ of $\Phi$ in the space $\fM(E_1,E_2)$ is, by definition, the norm of the linear transformer
$$
T\mapsto\int\limits_\X\int\limits_\Y\Phi(x,y)\,d E_1(x)T\,dE_2(y)
$$
on the class $\bS_1$.

If $\Phi\in\fM(E_1,E_2)$, one can define by duality double operator integrals of the form \rf{doi} for an arbitrary bounded linear operator $T$. However, we do not need this in this paper.

We are going to discuss briefly in \S~\ref{Haag} characterizations of Schur multipliers.

Birman and Solomyak proved in \cite{BS4} that if $f$ is a Lipschitz function and $A$ and $B$ are not necessarily bounded self-adjoint operators with $A-B\in\bS_2$, then 
\bay
\label{razS2}
f(A)-f(B)=\iint\limits_{\R\times\R}\frac{f(x)-f(y)}{x-y}\,dE_A(x)(A-B)\,dE_B(y).
\ey
Note that for an arbitrary Lischitz function $f$, the divided difference $\dg f$
is not always naturally defined on the diagonal. However, we can define $\dg f$
on the diagonal by an arbitrary bounded measurable function and the right-hand side
of \rf{razS2} does not depend on the values on the diagonal. It follows from \rf{razS2} that
$$
\|f(A)-f(B)\|_{\bS_2}\le\|f\|_{\Li}\|A-B\|_{\bS_2},
$$
where the Lipschitz (semi)norm $\|f\|_{\Li}$ of $f$ is, by definition,
$$
\sup\left\{\frac{|f(x)-f(y)|}{|x-y|}:~x,~y\in\R,~x\ne y\right\}.
$$

On the other hand, if $A$ and $B$ are not necessarily bounded self-adjoint operators with $A-B\in\bS_1$ and $f$ is an operator Lipschitz function, then
\bay
\label{razS1}
f(A)-f(B)=\iint\limits_{\R\times\R}\big(\dg f\big)(x,y)\,dE_A(x)(A-B)\,dE_B(y).
\ey
(see \cite{BS4}). Here the divided difference $\dg f$ is defined by \rf{razdraz}.
It follows from \rf{razS1} that
$$
\|f(A)-f(B)\|_{\bS_1}\le\|f\|_{\OL}\|A-B\|_{\bS_1}.
$$

\

\section{\bf Differentiability in the Hilbert--Schmidt norm}
\setcounter{equation}{0}
\label{HSdif}

\

Suppose that $A$ and $B$ are not necessarily bounded self-adjoint operators on Hilbert space such that $A-B\in\bS_2$. Consider the parametric family $A_t$, $0\le t\le1$, defined by $A_t=A+tK$, where $K\df B-A$. We need the following result of \cite{KPSS}, Theorem 7.18:

{\it Suppose that $f$ is a Lipschitz function on $\R$ that is differentiable at every point of $\R$. Then the function $s\mapsto f(A_s)-f(A)$ is differentiable on $[0,1]$ in the Hilbert--Schmidt norm and
$$
\frac{d}{ds}\big(f(A_s)-f(A)\big)\Big|_{s=t}
=\iint\limits_{\R\times\R}(\dg f)(x,y)\,dE_t(x)K\,dE_t(y).
$$
}

\

\section{\bf Schur multipliers and Haagerup tensor products}
\setcounter{equation}{0}
\label{Haag}

\

Let $(\X,E_1)$ and $(\Y,E_2)$ be spaces with spectral measures $E_1$ and $E_2$
on Hilbert space. There are several characterizations of the class $\fM(E_1,E_2)$ 
of Schur multipliers, see \cite{Pe2}, \cite{Pi}, \cite{AP}. We need the following characterization in terms of the Haagerup tensor product of $L^\be$ spaces:

{\it Let $\Phi$ be a measurable function on $\X\times\Y$. Then $\Phi\in\fM(E_1,E_2)$ if and only if $\Phi$ belongs to the Haagerup tensor product 
$L^\be(E_1)\!\otimes_{\rm h}\!L^\be(E_2)$, i.e., $\Phi$ admits a representaion
$$
\Phi(x,y)=\sum_n\f_n(x)\psi_n(y),
$$
where $\f_n\in L^\be(E_1)$, $\psi_n\in L^\be(E_2)$, and}
$$
\sum_n|\f_n|^2\in L^\be(E_1)\quad\mbox{and}\quad\sum_n|\psi_n|^2\in L^\be(E_2).
$$

Suppose now that $E_1$ and $E_2$ are Borel spectral measures on locally compact topological spaces $X$ and $Y$. In this case the following result holds (see \cite{AP}):

{\it Let $\Phi$ be a function on $\X\times\Y$ that is continuous in each variable.
Then \lb$\Phi\in\fM(E_1,E_2)$ if and only if it belongs to the Haagerup tensor product
$C_{\rm b}(\X)\!\otimes_{\rm h}\!C_{\rm b}(\Y)$ of the spaces of bounded continuous functions on $\X$ and $\Y$, i.e., $\Phi$ admits a representation 
$$
\Phi(x,y)=\sum_n\f_n(x)\psi_n(y),
$$
where $\f_n\in C_{\rm b}(\X)$, $\psi_n\in C_{\rm b}(\Y)$ and the functions
$$
\sum_n|\f_n|^2\quad\mbox{and}\quad\sum_n|\psi_n|^2
$$
are bounded.}

\

\section{\bf The trace of double operator integrals}
\setcounter{equation}{0}
\label{sleddoi}

\

Let $T$ be a trace class operator on Hilbert space and let
$E$ be a spectral measure on a $\s$-algebra of subsets of a set $\X$.
If $\Phi$ is a Schur multiplier, i.e., $\Phi\in\fM(E,E)$, then the double operator integral
$$
\iint\Phi(x,y)\,dE(x)T\,dE(y)
$$
belongs to $\bS_1$. Let us compute its trace. In \cite{BS4} the following trace formula
was found:
\bay
\label{BStf}
\trace\left(\iint\Phi(x,y)\,dE(x)T\,dE(y)\right)=\int\Phi(x,x)\,d\mu(x),
\ey
where $\mu$ is the complex measure defined by
$$
\mu(\D)=\trace\big(TE(\D)\big).
$$

The problem is how we can interpret the function $x\mapsto\Phi(x,x)$
for functions $\Phi$ in $\fM(E,E)$. In 
\cite{Pe5} the following justification of formula \rf{BStf} was given. We can define the trace ${\mathscr T}\Phi$ of a function $\Phi$ in $\fM(E,E)$  on the diagonal by the formula
$$
({\mathscr T}\Phi)(x)\df\sum_n\f_n(x)\psi_n(x),
$$
where 
\bay
\label{Htpre}
\Phi(x,y)=\sum_n\f_n(x)\psi_n(y)
\ey
is a representation of $\Phi$ as an element of the Haagerup tensor product 
$L^\be(E)\otimes_{\rm h}L^\be(E)$, i.e.,
\bay
\label{Htner}
\sum_m|\f_n|^2\in L^\be(E)\quad\mbox{and}\quad\sum_m|\psi_n|^2\in L^\be(E).
\ey
Clearly, the trace of $\Phi\in\fM(E,E)$ on the diagonal belongs to $L^\be(E)$.
Then formula \rf{BStf} holds if $\Phi(x,x)$ is understood as $({\mathscr T}\Phi)(x)$,
see \cite{Pe5}, \S~1.1.

Suppose now that $E$ is a Borel spectral measure on a locally compact topological space 
$\X$ and $\Phi$ is a function on $\X\times\X$ that is continuous in each variable.
As we have mentioned in \S~\ref{dois}, $\Phi$ admits a representation of the form
\rf{Htpre}, in which the functions $\f_n$ and $\psi_n$ satisfy \rf{Htner} and are continuous functions on $\X$. It is easy to see that in this case 
$({\mathscr T}\Phi)(x)=\Phi(x,x)$, $x\in\X$. In other words, the following theorem holds:

\begin{thm}
\label{sleddvoi}
Let $E$ be a spectral measure on a locally compact topological space $\X$. Suppose that 
$\Phi$ is a function of class $\fM(E,E)$. If $\Phi$ is continuous in each variable, then
formula {\em\rf{BStf}} holds for an arbitrary trace class operator $T$.
\end{thm}

\

\section{\bf The Lifshits--Krein trace formula for \\ arbitrary operator Lipschitz functions}
\setcounter{equation}{0}
\label{self}

\

Suppose that $A$ and $B$ are self-adjoint operators on Hilbert space such that
$A-B\in\bS_1$. Let $\bs\xi$ be the spectral shift function associated with the pair
$(A,B)$. As we have mentioned in \S~\ref{dois}, for an arbitrary operator Lipschitz function $f$
on $\R$, the operator $f(A)-f(B)$ belongs to trace class. The following theorem is the main result of the paper.

\begin{thm}
\label{glavn}
Let $f\in\OL(\R)$. Then
$$
\trace\big(f(A)-f(B)\big)=\int_\R f'(s)\bs{\xi}(s)\,ds.
$$
\end{thm}

To prove the theorem, we are going to use an approach of Birman and Solomyak in \cite{BSsl} to the Lifshits--Krein trace formula. In \cite{BSsl} they used their approach under more restrictive assumptions on $f$.

\medskip

\Pf Obviously, $f$ is a Lipschitz function. As we have mentioned in the introduction, $f$ is a differentiable function at every point of $\R$ (but not necessarily continuously differentiable!). 
Put $K\df B-A$.
Consider the parametric family $\{A_t\}_{0\le t\le1}$, $A_t\df A+tK$. Then $A_0=A$ and $A_1=B$. The operator $K$ obviously belongs to the Hilbert--Schmidt class $\bS_2$. As we have mentioned in \S~\ref{HSdif}, the function
\mbox{$t\mapsto f(A_t)-f(A)$} is differentiable in the Hilbert--Schmidt norm and
$$
Q_t\df
\frac{d}{ds}\big(f(A_s)-f(A)\big)\Big|_{s=t}=
\iint\limits_{\R\times\R}\frac{f(x)-f(y)}{x-y}\,dE_t(x)K\,dE_t(y)\in\bS_2,
$$
where $E_t$ is the spectral measure of $A_t$.

On the other hand, since the divided difference $\dg f$ is a Schur multiplier of of $\bS_1$ (see the introduction), it follows that
$$
Q_t\in\bS_1,\quad 0\le t\le1,\quad\mbox{and}\quad\sup_{t\in[0,1]}\|Q_t\|_{\bS_1}<\be.
$$

We have
$$
f(A)-f(B)=-\int_0^1Q_t\,dt,
$$
where the integral on the right is understood in the sense of Bochner in the space $\bS_1$. It follows that
$$
\trace\big(f(A)-f(B)\big)=-\int_0^1\trace(Q_t)\,dt.
$$

Since the function $f$ is differentiable everywhere, the divided difference $\dg f$ 
is continuous in each variable. By Theorem \ref{sleddvoi},
$$
\trace Q_t=\int_\R f'(x)\,d\nu_t(x),
$$
where the signed measure $\nu_t$ is defined by
$$
\nu_t(\D)\df \trace\big(E_t(\D)K)\quad\mbox{for a Borel subset}\quad\D\quad\mbox{of}\quad\R.
$$ 

We identify here the space of complex Borel measures on $\R$ with the dual space to the Banach space of continuous functions on $\R$ with zero limit at infinity. Then the function $t\mapsto \nu_t$ is continuous in the weak-$*$ topology on the space of complex Borel measures. Indeed, if $h$ is continuous on $\R$ and $\lim_{|x|\to\be}h(x)=0$, then
$$
\int h\,d\nu_t=\trace\big(h(A_t)K\big).
$$ 
The function $t\mapsto h(A_t)$ is a continuous function on $[0,1]$ in the operator norm; this follows from the fact that $h$ is an operator continuous function (see \cite{AP2}, \S~8). Thus 
the function $t\mapsto\trace\big(h(A_t)K\big)$ is continuous.

Therefore we can define the signed Borel measure $\nu$ on $\R$ by
$$
\nu=\int_0^1\nu_t\,dt.
$$
It follows that
$$
\trace\big(f(A)-f(B)\big)=-\int_\R f'\,d\nu.
$$

On the other hand, for smooth functions $g$ with compact support,
$$
\trace\big(g(A)-g(B)\big)=\int_\R g'\bs{\xi}\,d\m,
$$
where $\bs{\xi}$ is the spectral shift function. 

It follows that $\nu$ is absolutely continuous with respect to Lebesgue measure and
\lb$d\nu=-\bs{\xi}d\m$. $\bl$

\medskip

Theorem \ref{glavn} implies the following result:

\begin{thm}
Let $f$ be an operator Lipschitz function and let $A$ and $B$ be self-adjoint operators such that $A-B\in\bS_1$. Then the function
$$
t\mapsto\trace\big(f(A-tI)-f(B-tI)\big),\quad t\in\R,
$$
is continuous on $\R$.
\end{thm}

\Pf Consider the function $f_t$, $t\in\R$, defined by $f_t(x)\df f(x-t)$. Let $\bs{\xi}$
be the spectral shift function associated with $(A,B)$. It is easy to see that
$$
\trace\big(f(A+tI)-f(B+tI)\big)=\trace\big(f_t(A)-f_t(B)\big)
=\int_\R f'(x-t)\bs{\xi}(x)\,d\m(x),
$$
which depends on $t$ continuously. $\bl$

\

\

\

\footnotesize
\noindent
V.V. Peller \\
Department of Mathematics  \\
Michigan State University\\
East Lansing Michigan 48824\\

\end{document}